\documentclass[12pt]{article}
\usepackage{amssymb,amsmath,amsthm,amsfonts,enumerate, bbm, mathdots}

\usepackage{latexsym}
\usepackage{amscd}
\usepackage{stmaryrd}
\usepackage{color}
\usepackage[all]{xypic}
\usepackage{epsfig}
\usepackage{graphics}
\usepackage{ifthen}
\usepackage{varioref}
\usepackage{rotating}
\usepackage{extarrows}
\usepackage{cite}
\usepackage{mathrsfs}

\numberwithin{equation}{section}

\theoremstyle{plain}
\newtheorem{thm}{Theorem}[section]
\newtheorem{cor}[thm]{Corollary}
\newtheorem{lem}[thm]{Lemma}
\newtheorem{prop}[thm]{Proposition}
\newtheorem{rem}[thm]{Remark}
\newtheorem{Def}[thm]{Definition}
\newtheorem{eg}[thm]{Example}

\definecolor{darkgreen}{rgb}{0.0625,0.64,0.0625}
\usepackage[
            pdfstartview=FitH,
            CJKbookmarks=true,
            bookmarksnumbered=true,
            bookmarksopen=true,
            colorlinks,
            pdfborder=001,
            linkcolor=blue,
            anchorcolor=green,
            citecolor=red
            ]{hyperref}

\newfont{\scyr}{wncyr10 scaled 550}

\def\ker{\operatorname{Ker}}

\def\proof{\noindent {\bf Proof.\;}}
\def\Im{\operatorname{Im}}
\def\Hom{\operatorname{Hom}}

\def\id{\operatorname{id}}

\def\cha{\operatorname{char}}

\allowdisplaybreaks

\begin{document}

\title{Rota-Baxter systems of Hopf algebras and Hopf trusses}

\date{\small ~ \qquad\qquad School of Mathematical Sciences, Tongji University \newline No. 1239 Siping Road,
Shanghai 200092, China}

\author{Zhonghua Li\thanks{E-mail address: zhonghua\_li@tongji.edu.cn}~ and Shukun Wang\thanks{E-mail address: 2010165@tongji.edu.cn}}

\maketitle

\begin{abstract}
As Hopf truss analogues of Rota-Baxter Hopf algebras, the notion of Rota-Baxter systems of Hopf algebras is proposed. We study the relatiohship between Rota-Baxter systems of Hopf algebras and Rota-Baxter Hopf algebras, show that there is a Rota-Baxter system structure on the group algebra if the group has a Rota-Baxter system structure, investigate the descendent Hopf algebra of a Rota-Baxte system of Hopf algebras. Finally we study the local decomposition of the character group from a Rota-Baxter system of Hopf algebras to a commutative algebra.
\end{abstract}

{\small

}


\section{Introduction}\label{Sec:Intro}

\newcommand{\RNum}[1]{\uppercase\expandafter{romannumeral #1\relex}}

This paper arose from an attempt to understand the origin of Rota-Baxter operators on cocommutative Hopf algebras\cite{MG}. We investigate the Hopf truss\cite{TB2} ananlogues of Rota-Baxter Hopf algebras, study the relationships between Rota-Baxter Hopf algebras, Hopf braces\cite{IA} and Hopf trusses.

The (quantum) Yang-Baxter equation has been studied in mathematical physics since 1960s. In \cite{Drinfeld}, Drinfeld suggested to study the set-theoretical solution to the (quantum) Yang-Baxter equation. In 2007, the notion of braces \cite{Rump1} was introduced for abelian groups as tools to construct the set-solutions of the Yang-Baxter equation. In 2017, as noablian case of braces, skew braces were introduced to study non-involutive solutions in \cite{LG}.

 In 2017, Hopf braces were considered in \cite{IA}, as the quantum version of skew braces. Furthermore, it was shown that Hopf braces provide solutions to the (quantum) Yang-Baxter equation. In an attempt to understand the origins and the nature of the law binding two group operations together into skew braces, the notion of skew trusses was proposed in \cite{TB2} in 2019. As the linearization of skew trusses, Hopf trusses were also introduced in \cite{TB2}.

Rota-Baxter operators for associative algebras and Lie algebras arose in 1960s and they led to several interesting applications in renormalization of quantum field theory, solutions of classical Yang-Baxter equations, integrable system, operad theory and number theory, see for example \cite{CN}, \cite{SM} and see \cite{LG2008} for more details. In 2020, Rota-Baxter groups were introduced in \cite{LG1}, which are the group analogues of Rota-Baxter algebras and Lie algebras.  In 2021, Rota-Baxter Hopf algebras were proposed as the generalizations of Rota-Baxter groups and Lie algebras in \cite{MG}.

Recently, it was shown in \cite{VG} that every Rota-Baxter group is a skew brace and every skew brace can be embedded into a Rota-Baxter group. It was proven in \cite{YNL} that every Rota-Baxter Hopf algebra has the structure of Hopf braces. Conversely, it was shown in {\color{red} \cite{HHZ}} that a Hopf brace can be embedded into a Rota-Baxter Hopf algebra. In \cite{ZH}, the notion of Rota-Baxter systems of groups was proposed as the skew truss analogues of Rota-Baxter groups. Hence it is natural to consider the Hopf truss analogues of Rota-Baxter operators on cocommutative Hopf alegebras.	

Here is the outline of this paper. In Section \ref{2.0}, we mainly recall some basic definitions of Hopf braces and Rota-Baxter Hopf algebras. In Section \ref{4}, we define Rota-Baxter systems of Hopf algebras, invesgating the relationship between Rota-Baxter systems of Hopf algebras and Rota-Baxter Hopf algebras, showing there are descendent Hopf algebras induced by such algebraic systems. In Section \ref{5}, we derive the decomposition theorem of a character group from a Rota-Baxter system of Hopf algebras to a commutative algebra.

Throughout this paper, $\mathbb{F}$ is a field of characteristic $0.$ Without further mention, all vector spaces,
tensor products, homomorphisms, algebras, coalgebras, bialgebras, Lie algebras, Hopf algebras live over the field $\mathbb{F}$.
		
\section{Hopf braces and Rota-Baxter Hopf algebras}\label{2.0}

In this section, we recall the definitions of Hopf braces, Hopf trusses, Rota-Baxter Hopf algebras, Rota-Baxter systmes of groups and of Lie algebras.

\subsection{Hopf algebras}

We first list some basic definitions about Hopf algebras. For more details on Hopf algebras, we refer to \cite{SM,ME}.

For an algebra $A$ with the identity $1=1_A$, we will also denote the unit map $\mathbb{F}\to A; \alpha\mapsto \alpha 1$ by $1$. The multiplication map of $A$ is simply denoted by $\cdot$. Hence we may use the triple $(A,\cdot,1)$ to denote the algebra $A$.

A coalgebra is a vector space $C$ equipped with two linear maps, comultiplication $\Delta:C\to C\otimes C$ and counit $\epsilon:C\to \mathbb{F}$, such that
\begin{itemize}
  \item[(i)] $\Delta$ is coassociative: $(\Delta\otimes \id)\circ\Delta=(\id\otimes \Delta)\circ\Delta$,
  \item[(ii)] $\epsilon$ satisfies the counit property: $(\epsilon\otimes \id)\circ \Delta=\id=(\id\otimes \epsilon)\circ\Delta$.
\end{itemize}
The coalgebra $C$ is cocommutative if $\tau\circ\Delta=\Delta$, where $\tau:C\otimes C\to C\otimes C$ is the twist map defined by $a\otimes b\mapsto b\otimes a$.

For a coalgebra $C$, we will use sumless Sweedler
notation $\Delta(a)=a_{1}\otimes a_{2}$ with suppressed summation sign. When $\Delta$ is applied again to the left tensorand this would be written as $a_{1,1}\otimes a_{1,2}\otimes a_2$, while applying $\Delta$ to the right as $a_1\otimes a_{2,1}\otimes a_{2,2}$. Coassociativity means that these two expressions  are equal and hence it makes sense to write this element as $a_1\otimes a_2\otimes a_3$. Iterating this procedure gives
$$\Delta_{n-1}(a)=a_1\otimes a_2\otimes\cdots\otimes a_n,$$
where $\Delta_{n-1}(a)$ is the unique element obtained by applying coassociativity $(n-1)$ times. In this notation the counicity can be written as
$$\epsilon(a_1)a_2=a=a_1\epsilon(a_2),\quad\forall a\in C.$$

Let $(C,\Delta,\epsilon)$  and $(C',\Delta',\epsilon')$ be two coalgebras. A linear map $\phi:C\to C'$ is called a coalgebra homomorphism if for any $a\in C,$
\begin{align*}
\Delta'(\phi(a))=\phi(a_1)\otimes \phi(a_2),\quad \epsilon'(\phi(a))=\epsilon(a)
\end{align*}
and called a coalgebra anti-homomorphism if for any $a\in C$,
\begin{align*}
\Delta'(\phi(a))&=\phi(a_2)\otimes \phi(a_1),\quad \epsilon'(\phi(a))=\epsilon(a).
\end{align*}

Let $C$ be a coalgebra and $A$ be an algebra. The convolution product of $f,g$ in $\Hom(C,A)$ is defined to be the map $f\ast g$ given by
$$f\ast g(a)=f(a_1)g(a_2),\quad\forall a\in C.$$
Set $e=1_A\circ \epsilon_C$. Then the triple $(\Hom(C,A),\ast,e)$ is a unital algebra by \cite[Theorem 2.7]{LG2008}.

A bialgebra is a quintuple $(H,\cdot,1,\Delta,\epsilon)$ where $(H,\cdot,1)$ is an algebra and $(H,\Delta,\epsilon)$ is a coalgebra such that the multiplication map $\cdot$ and the unit map $1$ are coalgebra homomorphisms, which is equivalent to the comultiplication map $\Delta$ and the counit map $\epsilon$ are algebra homomorphisms.

A Hopf algebra is a bialgebra $(H,\cdot,1,\Delta,\epsilon)$ with a linear map $S:H\to H$ such that
$$a_1S(a_2)=\epsilon(a)1=S(a_1)a_2,\quad\forall a\in H.$$
Then map $S$ is called an antipode for $H$.

A map $\phi:H\to H'$ of bialgebras (Hopf algebras) $H,H'$ is called a bialgebra (Hopf algebra) homomorphism if it is both an algebra and a coalgebra homomorphism (and $\phi\circ S=S'\circ \phi$).

Let $(H,\cdot,1,\Delta,\epsilon,S)$ be a Hopf algebra. It is easy to see that if $H$ is either commutative or cocommutative, then $S^2=\id$. Furthermore, one can show that $S$ is an algebra anti-homomorphism and a coalgebra anti-homomorphism.

\subsection{Hopf braces and Hopf trusses}

The definition of Hopf braces was given in \cite{IA}.

\begin{Def}
 A Hopf brace over a coalgebra $(H,\Delta,\epsilon)$ consists of two Hopf algebra structures, denoted by $(H,\cdot,1,\Delta,\epsilon,S)$ and respectively by $(H,\circ,1_{\circ},\Delta,\epsilon,T)$, which satisfy the following compatibility
\begin{equation}\label{HB}
a\circ(bc)=(a_1\circ b)S(a_2)(a_3\circ c), \quad \forall a,b,c\in H.
\end{equation}
\end{Def}

As in \cite{IA}, we will denote $(H,\cdot,1,\Delta,\epsilon,S)$ by $H$, $(H,\circ,1_{\circ},\Delta,\epsilon,T)$ by $H_{\circ},$ and the Hopf brace by $(H,\cdot,\circ).$ Note that in any Hopf brace, $1_{\circ}=1.$

Let $(H,\cdot,\circ)$ and $(H',\cdot,\circ)$ be two Hopf braces. A linear map $f:(H,\cdot,\circ)\to (H',\cdot,\circ)$  is called a Hopf brace homomorphism if both $f:H\to H'$ and $f:H_{\circ}\to H_{\circ}'$ are Hopf algebra homomorphisms.

A Hopf brace $(H,\cdot,\circ)$ is called cocommutative if its underlying coalgebra $(H,\Delta,\epsilon)$ is cocommutative. It was proven in  \cite{IA} that the cocommutative Hopf braces provide solutions to the quantum Yang-Baxter equation.

As generalizations of skew trusses and Hopf braces, Hopf trusses were introduced in \cite{TB2}.

\begin{Def}
Let $(H,\cdot,1,\Delta,\epsilon,S)$ be a Hopf algebra. Let $\circ$ be a binary operation on $H$ making $(H,\circ,\Delta,\epsilon)$ a nonunital
bialgebra. We say that $(H,\cdot,\circ)$ is a Hopf truss if there is a coalgebra homomorphism $\sigma:H\to H$ such that
\begin{equation}\label{HTS}
a\circ(b c)=(a_1\circ b)S(\sigma(a_2))(a_3\circ c), \quad\forall a,b,c\in H.
\end{equation}
The map $\sigma$ is called the cocycle of the Hopf truss $(H,\cdot,\circ)$.
\end{Def}

The Hopf truss is a Hopf brace if $(H,\circ,\Delta,\epsilon)$ becomes a Hopf algebra and its cocycle $\sigma$ is the identity map.

\subsection{Rota-Baxter Hopf algebras}

Now we recall the notion of Rota-Baxter Hopf algebras introduced in \cite{MG}.

\begin{Def}
Let $(H,\cdot,1,\Delta,\epsilon,S)$ be a cocommutative Hopf algebra. A coalgebra homomorphism $B:H\to H$ is called a \textit{Rota-Baxter operator} on $H$ if for any $a,b\in H,$
\begin{equation}\label{HRB}
B(a)B(b)=B(a_1B(a_2)bS(B(a_3))).
\end{equation}
 Moreover, the pair $(H,B)$ is called a \textit{Rota-Baxter Hopf algebra}.
\end{Def}

For a Rota-Baxter Hopf algebra $(H,B)$, define an operation $\circ:H\otimes H \to H$ by
\begin{equation*}
	a\circ b=a_1B(a_2)bS(B(a_{3})), \quad\forall a,b\in H.
\end{equation*}
And define a map $T:H\to H$ by $T(a)=S(B(a_1))S(a_2)B(a_3)$ for any $a\in H.$  Then by \cite[Theorem 3]{MG}, $(H,\circ,1,\Delta,\epsilon, T)$ is also a cocommutative Hopf algebra, called the descendent Hopf algebra of $(H,B)$. It was shown in \cite[Theorem 2.13]{YNL} that $(H,\cdot,\circ)$ is a Hopf brace.

\subsection{Rota-Baxter systems of groups and of Lie algebras}

We recall the definitions of Rota-Baxter systems of groups and of Lie algebras from \cite{ZH}.

A Rota-Baxter system of groups is a group $G$ with two operators $B_1:G\to G$ and $B_2:G\to G,$ such that
\begin{align*}
B_1(a)B_1(b)&=B_1(B_1(a)bB_2(a)),\quad \forall a,b\in G,\\
B_2(b)B_2(a)&=B_2(B_1(a)bB_2(a)),\quad \forall a,b\in G.
\end{align*}
A Rota-Baxter system of Lie algebras is a triple $(\mathfrak{g},B_1,B_2)$ consists of a Lie algebra $\mathfrak{g}$ and two linear operators $B_1:\mathfrak{g}\to \mathfrak{g}$ and $B_2:\mathfrak{g}\to \mathfrak{g}$ such that
\begin{subequations}
	\begin{align*}
		&[B_1(a),B_1(b)]=B_1([B_1(a),B_1(b)]-[B_2(a),B_2(b)]),\quad\forall a,b\in \mathfrak{g},\\
		&[B_2(b),B_2(a)]=B_2([B_1(a),B_1(b)]-[B_2(a),B_2(b)]),\quad\forall a,b\in \mathfrak{g}.
	\end{align*}
\end{subequations}

\section{Rota-Baxter system of Hopf algebras}\label{4}

In this section, we define Rota-Baxter systems of Hopf algebras, as generalizations of Rota-Baxter Hopf algebras and Rota-Baxter systems of groups. We study the relationship between the Rota-Baxter systems of Hopf algebras and the Rota-Baxter systems of groups. Then as in the case for Rota-Baxter systems of groups, we study the descendent Hopf algebras of Rota-Baxter systems of Hopf algebras, which are induced from the descendent operations.

\subsection{Definition and examples}

\begin{Def}
Let $(H,\cdot,1,\Delta,\epsilon,S)$ be a cocommutative Hopf algebra with two coalgebra homomorphisms $B_1:H\to H$ and $B_2:H\to H$. The triple $(H,B_1,B_2)$ is called a Rota-Baxter system of Hopf algebras if $B_1(1)=B_2(1)=1$ and for any $a,b\in H$,
\begin{subequations}
\begin{align}
B_1(a)B_1(b)=B_1(B_1(a_1)bS(B_2(a_2))),\label{RB1}\\
B_2(a)B_2(b)=B_2(B_1(a_1)bS(B_2(a_2))).\label{RB2}
\end{align}
\end{subequations}
The operation $\circ:H\times H \to H$ given by
$$a\circ b=B_1(a_1)bS(B_2(a_2)), \quad\forall a,b\in H$$
is called the descendent operation of $(H,B_1,B_2).$
The map $\sigma:H\to H$ defined by
$$\sigma(a)=B_1(a_1)S(B_2(a_2)),\quad\forall a\in H$$
is called the cocycle of $(H,B_1,B_2)$.
\end{Def}

Note that $\sigma=B_1\ast (S\circ B_2)$, $a\circ 1=\sigma(a)$ and $1\circ a=a$ for any $a\in H$. Moreover, the following lemma holds.

\begin{lem}\label{IDM}
Let $(H,B_1,B_2)$ be a Rota-Baxter system of Hopf algebras with the cocycle $\sigma$. Then $B_1\circ \sigma=B_1$ and $B_2\circ \sigma=B_2$. Furthermore, $\sigma$ is idempotent and a coalgebra homomorphism.
\end{lem}

\proof
For any $a\in H$, by \eqref{RB1},
$$B_1(a)=B_1(a)B_1(1)=B_1(B_1(a_1)S(B_2(a_2)))=B_1(\sigma(a)).$$
Similarly, we get $B_2(\sigma(a))=B_2(a)$ from \eqref{RB2}.

Since $B_1$ is a coalgebra homomorphism and $S\circ B_2$ is a coalgebra anti-homomorphism, we have
\begin{align*}
\Delta(\sigma(a))&=\Delta(B_1(a_1)S(B_2(a_2)))=B_1(a_1)S(B_2(a_4))\otimes B_1(a_2)S(B_2(a_3)).
\end{align*}
Therefore
\begin{align*}
\Delta(\sigma(a))=B_1(a_1)S(B_2(a_2))\otimes B_1(a_3)S(B_2(a_4))=\sigma(a_1)\otimes \sigma(a_2).
\end{align*}
Then it is easy to see that $\sigma$ is a coalgebra homomorphism. Finally, it is immediate to see that $\sigma$ is idempotent.
\qed

In the sequel, we give some examples of Rota-Baxter systems of Hopf algebras.

\begin{eg}\label{egHopf}
If $(H,B)$ is a Rota-Baxter Hopf algebra, then we get a Rota-Baxter system of Hopf algebras $(H,B_1,B)$, where $B_1:H\to H$ is given by
$$B_1(a)=a_1B(a_2), \quad \forall a\in H.$$
\end{eg}

Conversely, we have the next proposition, which generalizes \cite[Corollary 3.14]{ZH}.

\begin{prop}
Let $(H,B_1,B_2)$ be a Rota-Baxter system of Hopf algebras with the cocycle $\sigma$.  If $\sigma$ is surjective, then $(H,B_2)$ and $(H,B_1\circ S)$ are Rota-Baxter Hopf algebras.
\end{prop}

\proof
By Lemma \ref{IDM}, $\sigma$ is just the identity map of $H$. Then we have $B_1\ast(S\circ B_2)=\id$. Since for any $a\in H$,
$$S(B_2(a_1))B_2(a_2)=S(B_2(a)_1)B_2(a)_2=\epsilon(B_2(a))1=\epsilon(a)1=e(a),$$
where $e=1\circ\epsilon$, we have $(S\circ B_2)\ast B_2=e$. Then $B_1=\id\ast B_2$. Therefore by \eqref{RB2} we find that $B_2$ is a Rota-Baxter operator on $H$. Moreover, we have
$$B_1\circ S=(\id\ast B_2)\circ S=S\ast(B_2\circ S).$$
Then by \cite[Proposition 1]{MG}, $(H,B_1\circ S)$ is also a Rota-Baxter Hopf algebra.
\qed

Now we introduce the notion of twisted Rota-Baxter operators on Hopf algebras, which generalizes the twisted Rota-Baxter operators on algebras and on groups \cite{ZH}.

\begin{Def}
Let $(H,\cdot,1,\Delta,\epsilon,S)$ be a cocommutative Hopf algebra with a coalgebra homomorphism $B:H\to H$ and a bialgebra homomorphism  $\phi:H\to H$. Then $B$ is called a $\phi$-twisted Rota-Baxter operator if
$$B(a)B(b)=B(B(a_1)bS(\phi(B(a_2)))),\quad \forall a,b\in H.$$
\end{Def}

\begin{prop}
 Let $(H,\cdot,1,\Delta,\epsilon,S)$ be a cocommutative Hopf algebra with a bialgebra homomorphism $\phi:H\to H$. For a $\phi$-twisted operator $B:H\to H$, define
 $L:H\to H$ by
 $$L(a)=\phi(B(a)),\quad\forall a\in H.$$ Then $(H,B,L)$ is a Rota-Baxter system of Hopf algebras.
\end{prop}

\proof
It is enough to show \eqref{RB2}. For any $a,b\in H,$ we have
\begin{align*}
L(a)L(b)&=\phi(B(a))\phi(B(b))=\phi(B(a)B(b))\\
&=\phi(B(B(a_1)bL(a_2)))=L(B(a_1)bL(a_2)),
\end{align*}
as desired.
\qed

\begin{eg}
Let $(H,\cdot,1,\Delta,\epsilon,S)$ be a cocommutative Hopf algebra. If  a coalgebra map $B:H\to H$ satisfies
$$B(a)B(b)=B(B(a)b),\quad\forall a,b\in H,$$
then it is easy to verify that $B$ is an $e$-twisted Rota-Baxter operator with $e=1\circ \epsilon$. Hence $(H,B,e)$ is a Rota-Baxter system of Hopf algebras by the above proposition.
\end{eg}

\subsection{Rota-Baxter systems and trusses}

We show in the next proposition that Rota-Baxter systems of Hopf algebras are Hopf trusses.

\begin{prop}\label{RBSH}
Let $(H,B_1,B_2)$ is a Rota-Baxter system of Hopf algebras. Then $(H,\cdot,\circ)$ is a Hopf truss with the cocycle $\sigma$.
\end{prop}

\proof
First we show that $\circ$ is a coalgebra homomorphism. In fact, for any $a,b\in H$, we have
\begin{align*}
\Delta(a\circ b)&=\Delta(B_1(a_1)bS(B_2(a_2)))=B_1(a_1)b_1S(B_2(a_4))\otimes B_1(a_2)b_2S(B_2(a_3))\\
&=B_1(a_1)b_1S(B_2(a_2))\otimes B_1(a_3)b_2S(B_2(a_4))=a_1\circ b_1\otimes a_2\circ b_2
\end{align*}
and
\begin{align*}
\epsilon(a\circ b)&=\epsilon(B_1(a_1)bS(B_2(a_2)))=\epsilon(B_1(a_1))\epsilon(b)\epsilon(S(B_2(a_2)))\\
&=\epsilon(a_1)\epsilon(b)\epsilon(a_2)=\epsilon(a)\epsilon(b).
\end{align*}

Next we check that $\circ$ satisfies the associative law. In fact, for any $a,b,c\in H,$ we have
\begin{align*}
a\circ(b\circ c)&=a\circ(B_1(b_1)cS(B_2(b_2)))\\
&=B_1(a_1)B_1(b_1)cS(B_2(b_2))S(B_2(a_2))\\
&=B_1(B_1(a_1)b_1S(B_2(a_2)))cS(B_2(B_1(a_3)b_2S(B_2(a_4))))\\
&=B_1(a_1\circ b_1)cS(B_2(a_2\circ b_2))=(a\circ b)\circ c.
\end{align*}

Finally, by Lemma \ref{IDM}, it remain to prove \eqref{HTS}. For any $a,b,c\in H,$ we have
\begin{align*}
(a_1\circ b)S(\sigma(a_2))(a_3\circ c)=&(B_1(a_1)bS(B_2(a_2)))S(\sigma(a_3))(B_1(a_4)cS(B_2(a_5)))\\=&(B_1(a_1)bS(B_2(a_2)))B_2(a_4)S(B_1(a_3))(B_1(a_5)cS(B_2(a_6)))\\
=&(B_1(a_1)bS(B_2(a_2)))B_2(a_3)S(B_1(a_4))(B_1(a_5)cS(B_2(a_6)))\\
=&B_1(a_1)bcS(B_2(a_2))=a\circ(bc),
\end{align*}
as required.
\qed

\subsection{Group-like and primitive elements}

Let $(H,B_1,B_2)$ be a Rota-Baxter system of Hopf algebras. Set
$$G(H)=\{a\in H\mid a\text{\;is group-like}\},\quad P(H)=\{a\in H\mid a\text{\;is primitive}\}.$$
It is well known  that $G(H)$ is a group and $P(H)$ forms a Lie algebra. As $B_1$ and $B_2$ are coalgebra homomorphisms, for any $a\in H$,
\begin{itemize}
 \item[(1)] if $a\in G(H)$, then both $B_1(a)$ and $B_2(a)$ are in $G(H)$;
 \item[(2)] if $a\in P(H)$, then both $B_1(a)$ and $B_2(a)$ are in $P(H)$.
\end{itemize}
Hence we have restriction maps $B_1|_{G(H)}$, $B_1|_{P(H)}$, $B_2|_{G(H)}$ and $B_2|_{P(H)}$. We show in the following proposition that $G(H)$ (resp. $P(H)$) naturally becomes Rota-Baxter system of groups (resp. Lie algebras).

\begin{prop}
Let $(H,B_1,B_2)$ be a Rota-Baxter system of Hopf algebras. Then $(G(H),B_1|_{G(H)},S\circ B_2|_{G(H)})$ $($resp. $(P(H),B_1|_{P(H)},-B_2|_{P(H)}))$ is a Rota-Baxter system of groups $($resp. Lie algebras$)$.
\end{prop}

\proof
It is immediate to see that $(G(H),B_1|_{G(H)},S\circ B_2|_{G(H)})$ is a Rota-Baxter system of groups. For any $a,b\in P(H),$ using Lemma \ref{IDM}, we have
\begin{align*}
	&B_1(a)B_1(b)=B_1(a)B_1(\sigma(b))=B_1(B_1(a_1)\sigma(b)S(B_2(a_2)))\\
=&B_1(B_1(a)\sigma(b)+\sigma(b)S(B_2(a)))=B_1(B_1(a)\sigma(b)-\sigma(b)B_2(a)).
\end{align*}
Therefore
\begin{align*}
	[B_1(a),B_1(b)]&=B_1(a)B_1(b)-B_1(b)B_1(a)\\
	&=B_1(B_1(a)\sigma(b)-\sigma(b)B_2(a)-B_1(b)\sigma(a)+\sigma(a)B_2(b)).	
\end{align*}
Note that for any  $a\in P(H)$, $\sigma(a)=B_1(a_1)S(B_2(a_2))=B_1(a)-B_2(a)$. Thus we have
\begin{align*}
	[B_1(a),B_1(b)]&=B_1(B_1(a)B_1(b)-B_1(b)B_1(a)+B_2(b)B_2(a)-B_2(a)B_2(b))\\
	&=B_1([B_1(a),B_1(b)]-[B_2(a),B_2(b)]).
\end{align*}
Similarly, one can verify that
$$[B_2(b),B_2(a)]=-B_2([B_1(a),B_1(b)]-[B_2(a),B_2(b)]).$$
Hence $(P(H),B_1|_{P(H)},-B_2|_{P(H)})$ is a Rota-Baxter system of Lie algebras.
\qed

Conversely, for Rota-Baxter systems of groups, we have  the following theorem, which generalizes \cite[Theorem 1]{MG}.

\begin{thm}
Let $G$ be a group with the identity $1_G$ and the inverse map $\pi:G\to G.$ If $(G,B_1,B_2)$ is a Rota-Baxter system of groups such that $B_1(1_G)=B_2(1_G)=1_G$, then $B_1$ and $\pi\circ B_2$ can be uniquely extended to operators $\widetilde{B_1}:\mathbb{F}[G]\to \mathbb{F}[G]$ and $\widetilde{B_2}:\mathbb{F}[G]\to \mathbb{F}[G]$ on its group algebra $\mathbb{F}[G]$ respectively, such that $(\mathbb{F}[G],\widetilde{B_1},\widetilde{B_2})$ is a Rota-Baxter system of Hopf algebras.
\end{thm}

\proof
We can uniquely extend $B_1$ and $B_2$ on $\mathbb{F}[G]$ by
\begin{align*}
\widetilde{B_1}\left(\sum \alpha_ig_i\right)=\sum \alpha_i B_1(g_i), \\
\widetilde{B_2}\left(\sum \alpha_ig_i\right)=\sum \alpha_i B_2(g_i)^{-1}
\end{align*}
respectively, where $\alpha_i\in \mathbb{F}$ and $g_i\in G.$ For any $g=\sum \alpha_i g_i\in \mathbb{F}[G]$ and $h\in G,$
we have
\begin{align*}
\widetilde{B_2}(g)\widetilde{B_2}(h)&=\sum \alpha_iB_2(g_i)^{-1}B_2(h)^{-1}=\sum \alpha_iB_2(B_1(g_i)hB_2(g_i))^{-1}\\
&=\sum \alpha_iB_2(B_1(g_i)hS(B_2(g_i)^{-1}))^{-1}=\widetilde{B_2}(\widetilde{B_1}(g_1)hS(\widetilde{B_2}(g_2))).
\end{align*}
As elements of $G$ form a linear basis of $\mathbb{F}[G],$  \eqref{RB2} holds for any $g,h\in \mathbb{F}[G].$ It is similar to show \eqref{RB1}.\qed

\begin{rem}
By \cite[Theorem 2]{MG}, every Rota-Baxter operator of weight $1$ on a Lie algebra $\mathfrak{g}$ can be uniquely extended to a Rota-Baxter operator on its universal enveloping algebra $\operatorname{U}(\mathfrak{g})$. Then for a Rota-Baxter system $\mathfrak{g}$ of Lie algebras, it is natural to consider the problem that whether the Rota-Baxter system structure on $\mathfrak{g}$ can be extended to its universal envoloping algebra $\operatorname{U}(\mathfrak{g})$ or not?
\end{rem}

\subsection{Descendent Hopf algebras}

Let $(H,B_1,B_2)$ be a Rota-Baxter system of Hopf algebras with cocycle $\sigma$. Set $H_1=\Im(\sigma)$.

\begin{lem}\label{COA}
\begin{description}
	\item[(1)] $H_1\circ H_1\subseteq H_1$ and $\Delta(H_1)\subseteq H_1\otimes H_1.$
	\item[(2)] Denote the restrictions of $\circ,\Delta,\epsilon$ on $H_1$ by $\circ_1,\Delta_1,\epsilon_1$ respectively, then $(H_1,\circ_1,1)$ is a unital algebra and $(H_1,\Delta_1,\epsilon_1)$ is a cocommutative coalgebra.
\end{description}
\end{lem}

\proof
For any $a,b\in H,$ $$\sigma(a)\circ \sigma(b)=a\circ1\circ b\circ 1=a\circ b\circ 1=\sigma(a\circ b)\in H_1.$$
By Lemma \ref{IDM}, $\sigma$ is a coalgebra homomorphism, which implies $\Delta(H_1)\subseteq H_1\otimes H_1$. As the coalgebra $(H,\Delta,\epsilon)$ is cocommutative,
$(H_1,\Delta_1,\epsilon_1)$ is cocommutative.
\qed

 The following theorem generalizes \cite[Theorem 3]{MG}.

\begin{thm}
Let $(H,B_1,B_2)$ be a Rota-Baxter system of Hopf algebras. Define $T:H_1\to H_1$ by
$$T(a)=S(B_1(a_1))B_2(a_2),\quad\forall a\in H_1.$$
Then $H_{B_1,B_2}=(H_1,\circ_1,1,\Delta_1,\epsilon_1,T)$ is a cocommutative Hopf algebra. It is called the descendent Hopf algebra of $(H,B_1,B_2)$.
\end{thm}

\proof
By Lemma \ref{COA} and the proof of Proposition \ref{RBSH}, $(H_1,\circ_1,1,\Delta_1,\epsilon_1)$ is a bialgebra. It remains to prove that $T$ is an antipode, that is, for any $a\in H,$
$$\sigma(a_1)\circ T(\sigma(a_2))=T(\sigma(a_1))\circ \sigma(a_2)=\epsilon(\sigma(a))1=\epsilon(a)1.$$

Since $\sigma$ is a coalgebra homomorphism, $B_1\circ \sigma=B_1$ and $B_2\circ \sigma=B_2,$ we get
\begin{align*}
\sigma(a_1)\circ T(\sigma(a_2))&=B_1(\sigma(a_1))S(B_1(\sigma(a_3)))B_2(\sigma(a_4))S(B_2(\sigma(a_2)))\\
&=B_1(a_1)S(B_1(a_2))B_2(a_3)S(B_2(a_4))\\
&=\epsilon(B_1(a_1))\epsilon(B_2(a_2))1=\epsilon(a)1.
\end{align*}
Then we have
\begin{align*}
B_1(T(\sigma(a_1)))B_1(a_2)&=S(B_1(a_1))B_1(\sigma(a_2))B_1(T(\sigma(a_3)))B_1(a_4)\\
&=S(B_1(a_1))B_1(\sigma(a_2)\circ T(\sigma(a_3)))B_1(a_4)\\
&=S(B_1(a_1))\epsilon(a_2)B_1(a_3)=\epsilon(a)1.
\end{align*}
It is similar to show that
$$B_2(T(\sigma(a_1)))B_2(a_2)=\epsilon(a)1.$$
Therefore
\begin{align*}
T(\sigma(a_1))\circ \sigma(a_2)&=B_1(T(\sigma(a_1)))\sigma(a_3)S(B_2(T(\sigma(a_2))))\\
&=B_1(T(\sigma(a_1)))\sigma(a_2)S(B_2(T(\sigma(a_3))))\\
&=B_1(T(\sigma(a_1)))B_1(a_2)S(B_2(a_3))S(B_2(T(\sigma(a_4))))\\
&=\epsilon(a)1,
\end{align*}
which finishes the proof.
\qed

\begin{cor}\label{IMB1}
$B_1$ and $B_2$ are Hopf algebra homomorphisms from $H_{B_1,B_2}$ to $H.$ Moreover, $\Im(B_1)$ and $\Im(B_2)$ are Hopf subalgebras of $H.$
\end{cor}

\proof
By \eqref{RB1} and \eqref{RB2}, $B_1$ and $B_2$ are bialgebra homomorphisms from $H_{B_1,B_2}$ to $H$. Then using \cite[Lemma 4.0.4]{ME}, we have $B_1$ and $B_2$ are  Hopf algebra homomorphisms form $H_{B_1,B_2}$ to $H$.
\qed

\subsection{Graph characterization}

Let $(H,\cdot,1,\Delta,\epsilon,S)$ be a cocommutative Hopf algebra. There is a bialgebra structure on $H^{\otimes 3}=H\otimes H\otimes H$. The multiplication $\cdot_{H^{\otimes 3}}$ is given by
$$(a\otimes b\otimes c)\cdot_{H^{\otimes 3}} (a'\otimes b'\otimes c')=a_1a'\otimes \epsilon(b)a_2b'S(c_2) \otimes c_1c',\quad\forall a,b,c,a',b',c'\in H.$$
The comultiplication $\Delta_{H^{\otimes 3}}$ is defined as
$$\Delta_{H^{\otimes 3}}(a\otimes b\otimes c)=(a_1\otimes b_1\otimes c_1)\otimes(a_2\otimes b_2\otimes c_2),\quad\forall a,b,c\in H.$$
Set $1_{H^{\otimes 3}}=1\otimes 1\otimes 1$ and $\epsilon_{H^{\otimes 3}}=\epsilon\otimes\epsilon\otimes \epsilon$. Then one can readily verify that $(H^{\otimes 3},\cdot_{H^{\otimes 3}},1_{H^{\otimes 3}},\Delta_{H^{\otimes 3}},\epsilon_{H^{\otimes 3}})$ is a bialgebra. we denote this bialgebra by $H^{\otimes 3}_{\bullet}$.

\begin{Def}
Let $(H,\cdot,1,\Delta,\epsilon,S)$ be a cocommutative Hopf algebra. Let $f_1:H\to H$ and $f_2:H\to H$ be two coalgebra homomorphisms. Define the graph of $(f_1,f_2)$ by
	$$\operatorname{Gr}_{f_1,f_2}(H)=\{f_1(a_1)\otimes a_2\otimes f_2(a_3)|a\in H\}.$$
\end{Def}

Then we have the following theorem.

\begin{thm}
Let $(H,\cdot,1,\Delta,\epsilon,S)$ be a cocommutative Hopf algebra. Let $B_1:H\to H$ and $B_2:H\to H$ be two coalgebra homomorphisms with $B_1(1)=B_2(1)=1$. Then $(H,B_1,B_2)$ is a Rota-Baxter system of Hopf algebras if and only if
$\operatorname{Gr}_{B_1,B_2}(H)$ is a subbialgebra of $H^{\otimes 3}_{\bullet}$.
\end{thm}

\proof
Assume that $(H,B_1,B_2)$ is a Rota-Baxter system of Hopf algebras. For any $a,b\in H,$ we have
\begin{align*}
	&(B_1(a_1)\otimes a_2\otimes B_2(a_3))\cdot_{H^{\otimes 3}}(B_1(b_1)\otimes b_2\otimes B_2(b_3))\\
	=&B_1(a_1)B_1(b_1)\otimes \epsilon(a_3)B_1(a_2)b_2S(B_2(a_5))\otimes B_2(a_4)B_2(b_3) \\
	=&B_1(a_1)B_1(b_1)\otimes B_1(a_2)b_2S(B_2(a_3))\otimes B_2(a_4)B_2(b_3)\\
	=&B_1(a_1\circ b_1)\otimes a_2\circ b_2\otimes B_2(a_3\circ b_3),
\end{align*}
which means that $\operatorname{Gr}_{B_1,B_2}(H)$ is a subalgebra of $H^{\otimes 3}_{\bullet}.$
And we have
\begin{align*}
&\Delta_{H^{\otimes 3}}(B_1(a_1)\otimes a_2\otimes B_2(a_3))\\
=& (B_1(a_1)\otimes a_3 \otimes B_2(a_5))\otimes (B_1(a_2)\otimes a_4\otimes B_2(a_6))\\
=&(B_1(a_1)\otimes a_2 \otimes B_2(a_3))\otimes (B_1(a_4)\otimes a_5\otimes B_2(a_6)).
\end{align*}
Hence $\operatorname{Gr}_{B_1,B_2}(H)$ is a subbialgebra of $H^{\otimes 3}_{\bullet}$.

Conversely, assume that $\operatorname{Gr}_{B_1,B_2}(H)$ is a subbialgebra of $H^{\otimes 3}_{\bullet}$. For any $a,b\in H,$ there is a $t\in H$ such that
$$(B_1(a_1)\otimes a_2\otimes B_2(a_3))\cdot_{H^{\otimes 3}}(B_1(b_1)\otimes b_2\otimes B_2(b_3))=B_1(t_1)\otimes t_2 \otimes B_2(t_3).$$
Since
\begin{align*}
&(B_1(a_1)\otimes a_2\otimes B_2(a_3))\cdot_{H^{\otimes 3}}(B_1(b_1)\otimes b_2\otimes B_2(b_3))\\
=&B_1(a_1)B_1(b_1)\otimes B_1(a_2)b_2S(B_2(a_3))\otimes B_2(a_4)B_2(b_3),
\end{align*}
we get
$$B_1(t_1)\otimes t_2 \otimes B_2(t_3)=B_1(a_1)B_1(b_1)\otimes B_1(a_2)b_2S(B_2(a_3))\otimes B_2(a_4)B_2(b_3).$$
Applying $\epsilon\otimes \id\otimes\epsilon$, we find
$$t=B_1(a_1)bS(B_2(a_2))=a\circ b.$$
Applying $\id\otimes\epsilon\otimes \epsilon$, we get $B_1(t)=B_1(a)B_1(b)$. Therefore we have $B_1(a)B_1(b)=B_1(a\circ b)$. Similarly, applying $\epsilon\otimes \epsilon\otimes \id$, we get $B_2(a)B_2(b)=B_2(a\circ b)$. Then $(H,B_1,B_2)$ is a Rota-Baxter system of Hopf algebras.
\qed

\section{Decomposition theorem for character groups}\label{5}

In this section, let $(H,B_1,B_2)$ be a Rota-Baxter system of Hopf algebras with the descendent Hopf algebra $H_{B_1,B_2}.$ Let $A$ be an unital commutative algebra.

Recall from \cite{LG2008} that an element $f\in \Hom(H,A)$ is called a character if $f$ is an algebra homomorphism. And by \cite[Proposition 2.11]{LG2008}, the set of characters, denoted by $\cha(H,A),$ is a group with respect to $\ast$. The identity element of $\cha(H,A)$ is $e=1_A\circ \epsilon$. Denote the identity element of $\cha(H_{B_1,B_2},A)$ by $e'$.

Define $\mathcal{B}_1:\cha(H,A)\to\cha(H_{B_1,B_2},A)$ by
$$\mathcal{B}_1(f)(a)=f(B_1(a)),\quad\forall f\in \cha(H,A),\quad\forall a\in H_{B_1,B_2}.$$
Similarly, define  $\mathcal{B}_2:\cha(H,A)\to\cha(H_{B_1,B_2},A)$ by
$$\mathcal{B}_2(f)(a)=f(B_2(a)),\quad\forall f\in \cha(H,A),\quad\forall a\in H_{B_1,B_2}.$$
Then it is easy to show that

\begin{lem}\label{Hom}
$\mathcal{B}_1$ and $\mathcal{B}_2$ are group homomorphisms.
\end{lem}

Moreover, we have the following lemma.

\begin{lem}
$\mathcal{B}_1(\ker(\mathcal{B}_2))$ is a normal subgroup of $\Im(\mathcal{B}_1)$ and $\mathcal{B}_2(\ker(\mathcal{B}_1))$ is a normal subgroup of $\Im(\mathcal{B}_2)$.
\end{lem}

\proof
By Lemma \ref{Hom}, $\mathcal{B}_1(\ker(\mathcal{B}_2))$ is a subgroup of $\Im(\mathcal{B}_1)$. For any $f\in \cha(H,A)$ and $g\in \ker(\mathcal{B}_2)$, we have
$$\mathcal{B}_2(f\ast g \ast f^{\ast -1})=\mathcal{B}_2(f)\ast\mathcal{B}_2(g)\ast\mathcal{B}_2(f^{\ast -1})=\mathcal{B}_2(f)\ast\mathcal{B}_2(f^{\ast -1})=\mathcal{B}_2(e)=e'.$$
Hence $f*g*f^{\ast -1}\in \ker(\mathcal{B}_2)$, which implies
$$\mathcal{B}_1(f)\ast\mathcal{B}_1(g)\ast\mathcal{B}_1(f^{\ast -1})=\mathcal{B}_1(f\ast g\ast f^{\ast -1})\in \mathcal{B}_1(\ker(\mathcal{B}_2)).$$
Therefore $\mathcal{B}_1(\ker(\mathcal{B}_2))$ is a normal subgroup of $\Im(\mathcal{B}_1)$. It is similar to show that  $\mathcal{B}_2(\ker(\mathcal{B}_1))$ is a normal subgroup of $\Im(\mathcal{B}_2)$.
\qed

Similar to \cite{LG1}, define $\Theta:\Im(\mathcal{B}_1)/\mathcal{B}_1(\ker(\mathcal{B}_2))\to \Im(\mathcal{B}_2)/\mathcal{B}_2(\ker(\mathcal{B}_1))$ by
$$\Theta(\overline{\mathcal{B}_1(f)})=\overline{\mathcal{B}_2(f)},\quad \forall f\in \cha(H,A).$$
Now we verify $\Theta$ is well-defined. In fact, for any $g\in \ker(\mathcal{B}_2)$ and $f\in \cha(H,A),$ we have
\begin{align*}
\Theta(\overline{\mathcal{B}_1(f)\ast\mathcal{B}_1(g))}&=\Theta(\overline{\mathcal{B}_1(f\ast g)})\\
&=\overline{\mathcal{B}_2(f\ast g)}=\overline{\mathcal{B}_2(f)}\ast\overline{\mathcal{B}_2(g)}\\
&=\overline{\mathcal{B}_2(f)}=\Theta(\overline{\mathcal{B}_2(f)}).
\end{align*}
The operator $\Theta$ is called the Cayley transform.

\begin{prop}\label{CL}
The Cayley transform $\Theta$ is a group isomorphism.
\end{prop}

\proof
It is easy to see that $\Theta$ is surjective. Hence it is enough to show that $\Theta$ is injective. For any $f\in \cha(H,A)$, if there is a $g\in \ker(\mathcal{B}_1)$ such that $\mathcal{B}_2(f)=\mathcal{B}_2(g)$, then we have $$\mathcal{B}_2(f\ast g^{\ast -1})=\mathcal{B}_2(f)\ast \mathcal{B}_2(g)^{\ast -1}=e'.$$
Therefore $\mathcal{B}_2(f\ast g^{\ast -1})\in \ker(\mathcal{B}_2)$ and
$$\mathcal{B}_1(f)=\mathcal{B}_1(f\ast g^{\ast -1})\in \mathcal{B}_1(\ker(\mathcal{B}_2)),$$
which proves that $\Theta$ is injective.
\qed

\begin{lem}\label{PSI}
Define $\Psi:\cha(H,A)\to \cha(H_{B_1,B_2},A)$ by
	$$\Psi(f)(a)=f(a),\quad\forall  f\in \cha(H,A),\quad \forall a\in H_{B_1,B_2}.$$
Then $\Psi$ is a group homomorphism.
\end{lem}

\proof
It is enough to show that $\Psi(f)\in \cha(H_{B_1,B_2},A)$ for any $f\in \cha(H,A)$. In fact, for any $a,b\in H$, by Lemma \ref{IDM} we have
\begin{align*}
\Psi(f)(\sigma(a)\circ \sigma(b))&=f(B_1(a_1)\sigma(b)S(B_2)(a_2))\\
&=f(B_1(a_1))f(\sigma(b))f(S(B_2(a_2))).
\end{align*}
Then by the communtativity of $A,$ we have
\begin{align*}
\Psi(f)(\sigma(a)\circ \sigma(b))&=f(B_1(a_1))f(S(B_2(a_2)))f(\sigma(b))\\
&=f(B_1(a_1)S(B_2(a_2)))f(\sigma(b))\\
&=\Psi(f)(\sigma(a))\Psi(f)(\sigma(b)).
\end{align*}
Hence $\Psi(f)\in \cha(H_{B_1,B_2},A)$.
\qed

Consider the group $(\Im(\mathcal{B}_1)\times \Im(\mathcal{B}_2),\cdot_D)$ where the product $\cdot_D$ is given by
$$(f_1,f_2)\cdot_D (g_1,g_2)=(f_1\ast f_2,g_1\ast g_2),\quad \forall f_1,f_2\in \Im(\mathcal{B}_1),\quad \forall g_1,g_2\in \Im(\mathcal{B}_2).$$
Let $G_{\mathcal{B}_1,\mathcal{B}_2}$ denote the subset
$$G_{\mathcal{B}_1,\mathcal{B}_2}=\{(f_1,f_2)\in \Im(\mathcal{B}_1)\times \Im(\mathcal{B}_2)|\Theta(\overline{f_1})=\overline{f_2}\}.$$

\begin{lem}\label{ImPSI}
For any $(f_1,f_2)\in G_{\mathcal{B}_1,\mathcal{B}_2},$ we have $f_1\ast f^{\ast -1}_2\in \Im(\Psi)$.
\end{lem}

\proof
There is a $f\in \cha(H,A)$ such that $\mathcal{B}_1(f)=f_1$. As $\overline{\Theta(f_1)}=\overline{f_2},$ there is a $g\in \ker(\mathcal{B}_1)$ such that $\mathcal{B}_2(f)\ast \mathcal{B}_2(g)=f_2$. Then by Lemma \ref{Hom} we have
$$f_1=\mathcal{B}_1(f)=\mathcal{B}_1(f\ast g),\quad f_2^{\ast -1}=\mathcal{B}_2(f\ast g)^{\ast -1}.$$
Therefore for any $a\in H_{B_1,B_2}$,
\begin{align*}
f_1\ast f^{\ast -1}_2(a)&=\mathcal{B}_1(f\ast g)(a_1)\mathcal{B}_2(f\ast g)^{\ast -1}(a_2)=f\ast g(B_1(a_1))f\ast g(S(B_2(a_2)))\\
	&=f\ast g(\sigma(a))=f\ast g(a).
\end{align*}
Hence $f_1\ast f_2^{\ast -1}=\Psi(f\ast g)\in \Im(\Psi)$.
\qed

Based on the above lemma, define $\Phi:G_{\mathcal{B}_1,\mathcal{B}_2}\to \Im(\Psi)$ by
$$\Phi(f_1,f_2)=f_1\ast f_2^{\ast -1},\quad\forall (f_1,f_2)\in G_{\mathcal{B}_1,\mathcal{B}_2}.$$
We have the following proposition.

\begin{prop}\label{ISO}
$G_{\mathcal{B}_1,\mathcal{B}_2}$ is a subgroup of $(\Im(\mathcal{B}_1)\times \Im(\mathcal{B}_2),\cdot_D)$ and $\Phi$ is a group isomorphism.
\end{prop}

\proof
By Proposition \ref{CL}, $G_{\mathcal{B}_1,\mathcal{B}_2}$ is a subgroup of $(\Im(\mathcal{B}_1)\times \Im(\mathcal{B}_2),\cdot_D)$.
For any $(f_1,f_2),(g_1,g_2)\in G_{\mathcal{B}_1,\mathcal{B}_2},$ we have
\begin{align*}
\Phi(f_1,f_2)\ast \Phi(g_1,g_2)&=(f_1\ast f_2^{\ast -1})\ast (g_1\ast g_2^{\ast -1})=(f_1*g_1)*(g_2^{*-1}*f_2^{*-1})\\
&=\Phi(f_1\ast g_1,f_2\ast g_2)=\Phi((f_1,f_2)\cdot_{D} (g_1,g_2)).
\end{align*}
Therefore $\Phi$ is a group homomorphism. Assume that $(f_1,f_2)\in \ker(\Phi)$. By the proof of Lemma \ref{ImPSI}, there is a $f\in \cha(H,A)$ such that $\mathcal{B}_1(f)=f_1$ and $\mathcal{B}_2(f)=f_2.$ Then for any $a\in H_{B_1,B_2},$ we have
\begin{align*}
f(a)&=f(\sigma(a))=f(B_1(a_1)S(B_2(a_2)))\\
&=f(B_1(a_1))f(S(B_2(a_2)))=\mathcal{B}_1(f)(a_1)\mathcal{B}_2(f)^{\ast -1}(a_2)=e'(a),
\end{align*}
which implies that $\Phi$ is injective. For any $g\in \cha(H,A),$ we have $(\mathcal{B}_1(g),\mathcal{B}_2(g))\in G_{B_1,B_2}$ and
$$\Psi(g)=\mathcal{B}_1(g)\ast \mathcal{B}_2(g)^{\ast -1}=\Phi(\mathcal{B}_1(g),\mathcal{B}_2(g)),$$
which proves $\Phi$ is surjective.
\qed

Finally by Proposition \ref{ISO}, we get the decomposition theorem of character groups from a Rota-Baxter system of Hopf algebras to a commutative algebra. Note that the decomposition theorem for Rota-Baxter groups was obtained in \cite[Theorem 3.5]{LG1}.

\begin{thm}\label{DEM}
For any $f\in \Im(\Psi)$, there is a unique decomposition $f=f_1\ast f_2^{\ast -1}$  with $(f_1,f_2)\in G_{\mathcal{B}_1,\mathcal{B}_2}$.
\end{thm}

\begin{cor}
Let $(H,B)$ be a Rota-Baxter Hopf algebra. Then for any $f\in \cha(H,A)$, there is a unique decomposition $f=f_1\ast f_2$ with $(f_1,f_2)\in G_{\mathcal{B}_1,\mathcal{B}_2}$, where the Rota-Baxter system $(H,B_1,B_2)$ of Hopf algebras induced from $(H,B)$ is given in Example \ref{egHopf}.
\end{cor}

\proof
For the Rota-Baxter system of Hopf algebras $(H,B_1,B)$ given in Example \ref{egHopf}, it is easy to see $\sigma=\id_{H}$. Hence $H_{B_1,B}=H_{\circ}$ and it follows from Lemma \ref{PSI} that $\Psi$ is the identity map and $\cha(H,A)=\cha(H_{B_1,B},A).$ Then by Theorem \ref{DEM} we prove the assertion.
\qed

\end{document}